\documentclass[a4paper, 12pt]{article}

\def\margin{2cm}
\usepackage[left=\margin,right=\margin,top=\margin,bottom=\margin]{geometry}

\usepackage{authblk}

\usepackage{amsmath, amssymb, amsthm}
\usepackage{xcolor}
\usepackage{url}

\title{Upper bound on cubicity in terms of boxicity\\ for graphs of low chromatic number}

\author[1]{L.~Sunil~Chandran}
\author[2]{Rogers~Mathew\footnote{Supported by VATAT Post-doctoral Fellowship, Council of Higher Education, Israel.}}
\author[2]{Deepak~Rajendraprasad$^*$}
\affil[1]{
	Department of Computer Science and Automation, \authorcr 
	Indian Institute of Science, Bangalore, India - 560012. \authorcr
	\texttt{sunil@csa.iisc.ernet.in}
}
\affil[2]{
	The Caesarea Rothschild Institute, Department of Computer Science, \authorcr
	University of Haifa, 31095, Haifa, Israel. \authorcr
	\texttt{(rogersmathew,deepakmail)@gmail.com}
}

\theoremstyle{plain}
\newtheorem{thm}{Theorem}
\newtheorem{lem}[thm]{Lemma}
\newtheorem{cor}[thm]{Corollary}
\newtheorem{obs}[thm]{Observation}
\newtheorem{claim}[thm]{Claim}

\theoremstyle{definition}

\theoremstyle{remark}
\newtheorem*{remark}{Remark}

\newtheoremstyle{cases}
  {}
  {}
  {}
  {}
  {}
  {\newline}
  {0.5em}
   {{\itshape \thmname{#1}} \thmnumber{#2} ({\itshape\thmnote{#3}}).\medskip}
\theoremstyle{cases}
\newtheorem{case}{Case}

\newcommand{\boxi}{\ensuremath{\mathrm{box}}}
\newcommand{\cubi}{\ensuremath{\mathrm{cub}}}
\newcommand{\rogers}[1]{{\color{black} #1}}

\begin {document}

\pagestyle{plain}
\pagenumbering{arabic}
\maketitle

\begin {abstract}
 The boxicity (respectively cubicity) of a graph $G$ is the minimum non-negative integer $k$,
 such that $G$ can be represented as an  intersection graph of axis-parallel
 $k$-dimensional boxes (respectively $k$-dimensional  unit cubes) and is denoted by $\boxi(G)$ 
 (respectively  $\cubi(G)$). It was shown by Adiga and Chandran 
 (Journal of Graph Theory,
  65(4), 2010) that for any graph $G$, $\cubi(G) \le$ box$ (G) \left \lceil
   \log_2 \alpha \right \rceil$, where $\alpha = \alpha(G)$ is the 
  cardinality of the  maximum independent set in $G$.   
 In this note we show that $\cubi(G) \le 2 \left \lceil \log_2  \chi(G) \right
  \rceil \boxi(G) +   \chi(G)  \left \lceil \log_2 \alpha(G) \right \rceil $.
  In general, this result can provide a  much better upper bound than that of
  Adiga and Chandran for graph classes with bounded
  chromatic number. For example, for bipartite graphs we get,
  $\cubi(G) \le 2  (\boxi(G) +    \left \lceil \log_2 \alpha(G) \right \rceil )$. 

  Moreover we show that for every positive integer $k$, there exist graphs
  with chromatic number $k$, such that for every $\epsilon > 0$, 
  the value given by our upper bound is at most $(1+\epsilon)$ times 
  their cubicity. Thus, our upper bound is almost tight.

\end {abstract}

\section{Introduction}

An \emph{axis-parallel $k$-dimensional box}, or \emph{$k$-box} in short,
is the Cartesian product $R_1\times R_2\times\cdots\times R_k$ where each
$R_i$ is an interval of the form $[a_i,b_i]$ on the real line. 
A $k$-unit-cube (or $k$-cube in short)  is a  $k$-box where each $R_i$ is an interval of the form $[a_i,a_i +1]$.
A graph $G(V,E)$ is said to be an \emph{intersection graph} of $k$-boxes (respectively $k$-cubes) 
if there is a mapping $f$ that maps the vertices of $G$ to $k$-boxes (respectively $k$-cubes) 
such that for any two vertices $u,v\in V$, $(u,v)\in E(G)\Leftrightarrow
f(u)\cap f(v)\not=\emptyset$. Then  $f$ is called a \emph{$k$-box
representation} of $G$ (respectively $k$-cube representation). 
The \emph{boxicity} (respectively \emph {cubicity} ) of a graph $G$, denoted by $\boxi(G)$ (respectively  $\cubi(G)$)
is the minimum non-negative integer $k$ such that $G$ has a $k$-box representation
(respectively $k$-cube representation). Only complete graphs have boxicity (cubicity)  0. 
The class of graphs with boxicity at most 1 is the class of interval 
graphs, and the class of graphs  with cubicity at most 1  is the class of unit interval graphs.

In the rest of the paper we will always use $n$ to denote the number of
vertices of the graph being discussed. Logarithms will be to the base 2,
unless otherwise specified.

Let $H_1$ and $H_2$ be two graphs such that $V(H_1) = V(H_2) = V(G)$
and $E(G) = E(H_1) \cap E(H_2)$.  Then we write $G= H_1 \cap H_2$.
The following observation was made by F. S. Roberts \cite {Rob1}.

\begin {lem}[\cite{Rob1}]
\label {Robertslem} 
Boxicity of a non-complete graph $G$  is the minimum positive  integer $k$ such that there exists
interval graphs $I_1,\ldots, I_k$  such that
$G = I_1 \cap \cdots \cap I_k$.  Cubicity of a non-complete  graph $G$ is the 
minimum positive  integer $k$ such that there exists $k$ unit interval graphs
$U_1,\ldots, U_k$ such that 
$G = U_1 \cap \ldots \cap U_k$. 
\end {lem}

\subsection{Brief literature survey}

The concepts of boxicity and cubicity were introduced by F. S. Roberts \cite {Rob1} in 1968 for studying some problems in ecology.
 The computational complexity of
finding the boxicity of a graph was studied by \cite {Yan1,Coz1}. It is known that
it is NP-hard to decide whether the boxicity of a graph is at most 2 \cite {Krat1}.
Recently it was shown that it is hard to even approximate within 
$n^{1-\epsilon}$, for any $\epsilon > 0$ unless NP=ZPP  \cite {Chalermsook2013}. The best known approximation
factor for boxicity is  $O \left (\frac {n \sqrt  {\log \log n}  } {\sqrt  {\log n} } \right )$ 
and that for cubicity is $O \left (\frac {n (\log \log n )^{3/2} } {\sqrt { \log n } } \right  )$ \cite {AdigaIPEC}.

Adiga et al. \cite {AdigaCOCOON} showed that boxicity is closely  related to
the well-studied concept of  
partial order dimension. Chandran and Sivadasan
 \cite {Chandran2007} proved that for any graph $G$, $\boxi(G) \le \textnormal{treewidth}(G) + 2$. It is also known to be related
to graph minors: Esperet and Joret \cite {EsperetJoret13} 
 showed that  for any graph $G$, \rogers{$\boxi(G) \in O(\eta^4 \log^2 \eta)$}, where
$\eta $ is the number of vertices in the largest clique minor of $G$.
Chatterjee and Ghosh \cite {ChatterjeeGhosh}  related boxicity with  Ferrer's dimension.
The upper bound for $\boxi(G)$  in terms of maximum degree of $G$ (denoted by
$\Delta(G)$ ) was studied in \cite {ChandranFrancisSivadasan2008,Esperet09}. The current best known upper bound
for $\boxi(G)$ in terms of $\Delta(G) =\Delta$, is $O(\Delta \log^2 \Delta)$
\cite {AdigaCOCOON}, which follows from a corresponding upper bound for partial
order dimension in terms of the maximum degree of the comparability 
graph of the partial order \cite {FuerediKahn}. 

Boxicity of outerplanar graphs is known to be at most 2 \cite {Sch1}. 
Thomassen \cite {Thom1} showed that the boxicity of planar graphs is at most 
3. Recently Felsner and Francis \cite {Mathew2011} gave a different proof for the above
theorem. More  proofs of this theorem  are given in \cite {SueWhitesideetc}. 
Hartman et al.
showed that the boxicity of bipartite planar graphs
 is at most 2 \cite {Irith1991}.
Lower bounds for boxicity was studied in \cite {AdigaLB}. 

Chandran et al. \cite {ChandranFrancisSivadasan2013} showed that for any graph $G$, $\cubi(G) \le
b + 1$, where $b$ denotes the bandwidth of $G$. They also showed
that $\cubi(G) = O(\Delta \log b)$. Chandran and Mathew \cite {ChandranMathew2014} showed that
for any graph $G$, $\cubi(G) \le (k+2) \left \lceil 2e \log n \right \rceil$,
where $k$ is the degeneracy of $G$. Chandran and Sivadasan \cite {ChandranSivadasan2008}
showed that for $d$-dimensional hypercubes $H_d$, $\cubi(H_d) =
\theta \left ( \frac {d} {\log d } \right ) $. Adiga and Chandran
\cite {Adiga10} showed that for an interval graph $G$,
$\left \lceil \log \psi \right \rceil  \le \cubi(G) \le \left \lceil \log \psi \right \rceil + 2$,
where $\phi$ denotes the number of leaves in the largest induced star in $G$.  


\subsection{Cubicity vs. boxicity}
 
Clearly for any graph $G$, $ \boxi(G) \le \cubi(G)$. 
Then the following question becomes relevant: Does there exists a 
function $g$  such that $\cubi(G) \le g(\boxi (G))$? It is easy to see that
the answer is negative: Consider a star graph on $n+1$ vertices.
Its  cubicity is $\left \lceil \log n \right \rceil$ \cite {Rob1}  
whereas its
boxicity is 1, since a star graph is an interval graph.  
Chandran and Mathew \cite {Ashik2009} showed that for any graph $G$, 
$\cubi(G) \le   \boxi (G) \left \lceil \log n \right \rceil $, where $n$ is the
number of vertices.  
Adiga and Chandran \cite {Adiga10} improved this result by showing that
we can use the cardinality of the maximum independent set in $G$, denoted
by $\alpha(G)$,  in the place of $n$.

\begin {lem}[\cite{Adiga10}]
\label {Adigaslem}
For any graph $G$, $\cubi(G) \le \boxi(G) \left \lceil \log \alpha(G) 
\right \rceil$. 
\end {lem}

\begin{remark}
We demonstrate next that the bound in the above lemma is tight, i.e., given any two positive integers $b$ and $\alpha$, we show that there exists a graph $G$ with $\boxi(G) = b$, $\alpha(G) = \alpha$ and $\cubi(G) = b \lceil \log \alpha \rceil$. It is known that a complete $p$-partite graph, $p \geq 2$, with $n_i, i \in [p]$, vertices in each part has boxicity $p$ and cubicity $\sum_{i=1}^p \lceil \log n_i \rceil$ \cite{Rob1}. Hence, given positive integers $b$ and $\alpha$, if $b \geq 2$, the complete $b$-partite graph with $\alpha$ vertices in each part will serve our purpose. If $b = 1$, then a star graph on $\alpha+1$ vertices gives the required graph. 
\end{remark}

We observe that, in a loose sense, the two terms in the upper bound on $\cubi(G)$ 
given in Lemma \ref {Adigaslem}, 
namely $\left \lceil \log \alpha(G) \right \rceil$ and $\boxi(G)$,  
by themselves contribute in keeping the cubicity of a graph
high. Clearly $\boxi(G)$ is a lower bound for   $\cubi(G)$ since 
cubes are specialised boxes.
The other term  $\left \lceil \log \alpha(G) \right \rceil$ can make $\cubi(G)$ high  
due to a geometric reason, captured in the so-called
`volume argument',  which we reproduce here (also see \cite {ChandranOrioloMannino}): 
Let $\cubi(G) = k$. 
 Then the  vertices in the biggest
independent set should correspond to  pairwise non-intersecting $k$-dimensional
unit cubes in the
cube representation of $G$.  Thus if we consider the minimal bounding box
for the cube representation, that bounding box should have a volume of
at least $\alpha(G)$ units.  On the other hand, the width of this
bounding box on any of the dimensions (i.e. the distance between the
extreme points of the projection of the cube representation on the 
corresponding axis)  can be at most $d+1$ where $d$ is the diameter of $G$.
(To see this note that each vertex correspond to a unit length interval in the
projection and the pair of vertices containing the two extreme points in
their respective intervals, are at a distance of at most $d$ in the graph, and
thus it follows that the geometric distance between the two extreme points 
is at most $d+1$.)
From this it is clear that the volume of the minimal bounding box is at most
$(d+1)^k$. It follows that $(d+1)^k \ge  \alpha(G)$. From this
we get $\cubi(G) = k \ge \left \lceil \log_{d+1}  \alpha(G) \right \rceil$.
Thus we have $ M= \max ( \boxi(G),  \left \lceil \log_{d+1}  \alpha(G) \right \rceil ) \le \cubi(G)$.
It is natural to investigate whether there exists a function $g$ such that
$\cubi(G) \le g(M)$. But the answer is negative since we can increase the diameter
of a graph unboundedly without affecting its cubicity. For example, if $G$ is the graph obtained by identifying one end point of a path on $2n+1$ vertices with the leaf of a star graph on $n+1$ vertices, it is easy to check that $\boxi(G) = 1$, $\alpha(G) = 2n$, diameter of $G$ $d = 2n + 2$ and hence 
$M = \max \{\boxi(G), \lceil \log \alpha(G) / \log (d+1) \rceil \} = 1$,
whereas $\cubi(G) = \lceil \log n \rceil$ which is far higher.

In this paper we ask a simpler question:
Let $\bar{M} = \max ( \boxi(G), \left \lceil \log \alpha(G) \right \rceil)$. 
Lemma \ref{Adigaslem} tells us that $\cubi(G) \in O(\bar{M}^2)$ and the 
remark after the lemma indicates that we cannot have anything better in general
(choosing $\alpha = 2^b$ there illustrates the point). But can we show that $\cubi(G) \in O(\bar{M})$ for some restricted graph classes? 
In this paper we show that if we restrict  ourselves to classes of graphs
whose chromatic number is bounded above by a constant, such a result
can indeed be proved. 
In fact our main theorem is  a general upper bound for cubicity in terms of boxicity,
the independence number and chromatic number:

\begin {thm}
\label {mainthm}
Let $G$ be a graph with chromatic number $\chi(G)$ and the cardinality of
the maximum independent set $\alpha(G)$. Then 
$\cubi(G) \le 2 \left \lceil \log \chi(G)  \right \rceil \boxi(G) + 
   \chi(G)  \left \lceil \log \alpha(G)  \right \rceil $.
\end {thm}

 For graphs of low chromatic number, this result can  be in general  far better than that of Adiga et al. \cite{Adiga10}.
 The most interesting case is that of bipartite graphs:

\begin {cor}

\label {bipcorollary}
 For a bipartite graph $G$, $\cubi(G) \le 2 (\boxi(G) +  \left \lceil \log \alpha(G) \right \rceil)$.

\end {cor}

\begin{remark}
The reader may naturally wonder whether 
chromatic number is an upper bound for the boxicity of a graph or not,
in which case, Theorem \ref {mainthm} cannot be an improvement over
Lemma \ref {Adigaslem}. 
On the contrary there are several graphs with boxicity greater than
chromatic number. In fact it looks most graphs are like that: 
In \cite {AdigaLB} it is shown that almost all balanced
bipartite graphs (on $2n$ vertices)  have boxicity $\Omega(n)$.
The proof can be modified to show that almost all bipartite graphs 
with $n$ vertices on one side and $m$ vertices on the other, have
boxicity $\Omega (\min (n,m))$. 
\end{remark}

\subsection {Preliminaries}

A graph $G$ is a co-bipartite graph if the complement of it, $\overline G$
is a bipartite graph. Thus $G$ is a co-bipartite graph if and only if 
the vertex set $V(G)$ can be partitioned into two sets $A$ and $B$ such
that $A$ and $B$ both induce cliques in $G$. It is clear that
for a co-bipartite graph $G$, $\alpha(G) \le 2$. Applying Lemma
 \ref {Adigaslem}
we can infer the following:

\begin {lem}
\label {boxicitycubicityofcobipartitelem}
For a co-bipartite graph $G$, $\cubi(G) = \boxi(G)$. 
\end {lem}

Let $G$ be a graph and let its vertex set $V(G)$ be partitioned into 
$A$ and $B$. Now construct a graph $H$, with $V(H) = V(G)$
and $E(H) = E(G) \cup \{ (u,v) : u,v \in A \} \cup \{ (u,v): u,v \in B \}$.
 Note that $H$ is obtained from $G$ by adding more edges so that $A$ as well as $B$
induce cliques in $H$ and the edges across $A$ and $B$ are as in $G$. 
The following observation is from \cite {Jas1}.

\begin {lem}

\label {cobipartizinglem}

$\boxi(H) \le 2 \boxi(G)$ 

\end {lem}

\section {Proof of Theorem \ref {mainthm} }

Let $G$ be a graph with $\boxi(G) = b$.  Consider a proper coloring
of $G$ using $\chi(G)= \chi$ colors: Let $C_0,\ldots,C_{\chi-1}$ be
the color classes with respect to this coloring.  First we define 
$\lceil \log \chi \rceil$  bipartitions of $V(G)$ by the following rule:
The bipartition $(A_i,B_i)$ for $1 \le i \le \lceil \log \chi \rceil$ is
obtained by setting $A_i$ as the union of all the color classes $C_j$
such that  the $i$-th bit in the binary representation of $j$ is 1. 
(Here we consider binary representation of numbers in $\{0,1,\ldots,\chi-1\}$ 
using $\lceil \log \chi \rceil$  bits).  $B_i$ is defined as the
union of the remaining color classes.  Define $H_i$ to be the co-bipartite
graph obtained by defining  $V(H_i) = V(G)$  and
the edge set $E(H_i) = E(G) \cup \{ (u,v) : u,v \in A_i \} \cup
  \{ (u,v) : u,v \in B_i \}$. That is, $H_i$ is obtained by  adding  edges to $G$ such that
$A_i$ and $B_i$ induce  cliques in the resulting graph, and the edges
across $A_i$ and $B_i$ are as in $G$.  Since $H_i$ is a co-bipartite
graph, by Lemma \ref {boxicitycubicityofcobipartitelem},
 cub$(H_i)=\boxi(H_i)$. By Lemma \ref {cobipartizinglem},
the box$(H_i) \le 2 \boxi(G) = 2b$.  Thus cub$(H_i) \le 2b$. Therefore
by Lemma \ref {Robertslem} there exist $2b$ unit interval graphs, say
$U_i^1,\ldots,U_i^{2b}$ such that 

\begin {eqnarray}
\label  {cubicityofcobipartitegraph}
   U_i^1 \cap \cdots \cap U_i^{2b} = H_i
\end {eqnarray} 

\begin {obs}
\label {obs1}
 For $1 \le i \le \left \lceil \log \chi(G) \right \rceil $ and $1 \le j \le 2b$, $U_i^j$ is a super graph of $G$.
\end {obs}

\begin {proof}
 { From equation \ref {cubicityofcobipartitegraph} it is clear that $U_i^j$ is 
a super graph of $H_i$ which in turn is a super graph of $G$. }
\end {proof} 

Also define for each \rogers{color class $C_i$}, $0 \le i \le \chi-1$, 
$t_i = \lceil \log {|C_i|} \rceil$ unit interval graphs, $W_i^1,\ldots,W_i^{t_i}$
 in the following way: 
First number the vertices of $C_i$ from $0$ to $|C_i|-1$.
Let $n_i(u)$ be the number given to a vertex $u \in C_i$ as per this
numbering.
Now to define the unit interval graph $W_i^j$, for $0 \le i \le \chi-1$
and $1 \le j \le t_i$, associate
with each vertex $v \in V(G)$   an interval $f_i^j(v)$  as follows:
 
For each vertex $u \in V-C_i,  \ \ f_i^j(v) = [1,2]$ \\
 
For each $u \in C_i,  f_i^j(u) = [0,1]$   if the $j$-th bit
in the binary representation of $n_i(u)$ is $1$, else 
$f_i^j(u) = [2,3]$. Define $W_i^j$ to be the corresponding
unit interval graph. 

\begin {obs}
\label {obs2}
 For $0 \le i \le \chi - 1$, and $1 \le j \le t_i$, $W_i^j$ is a super graph of $G$
\end {obs}

\begin {proof}
 If $(u,v) \in E(G)$ then $u$ and $v$ do not belong to the same color class.
Let $u \in C_a$ and $v \in C_b$, where $a \ne b$.
Clearly one of $a,b \ne i$: Without loss of generality, let $b \ne i$. Thus
$f_i^j(v) = [1,2]$.
If $a \ne i$, then $ f_i^j(u) =  [1,2]$ and
if  $a=i$, then $f_i^j(u)$ equals either $[0,1]$ or $[2,3]$. In all cases,
$u$ is adjacent to $v$ in $W_i^j$. It follows that $W_i^j$ is a super graph
of $G$.   
\end {proof}

\begin {claim}
\label {mainclaim}
 $$ \bigcap_{1 \le i \le \left \lceil \log \chi \right \rceil ; 1 \le j \le 2b} U_i^j  \bigcap_{0 \le i \le \chi -1; 1 \le j \le t_i} W_i^j   = G $$ .  
\end {claim} 

In view of observations \ref {obs1} and \ref {obs2},  to prove the above claim it is sufficient to 
show that if $(u,v) \notin E(G)$ then there exists one unit
interval graph  $I \in \{ U_i^j :   1 \le i \le \left \lceil \log \chi \right \rceil ; 1 \le j \le 2b  \}  \cup  \{ W_i^j : 0 \le i \le \chi -1; 1 \le j \le t_i\}$ such that  $u$ is not adjacent to $v$ in it. 
We consider two cases:

\begin{case}[$u, v \in C_i$, for some $i \in \{0, \chi-1\}$] 
Recall that we had numbered the vertices of
$i$th color class from $0$ to $|C_i|-1$ to  define the interval
graphs $W_i^j$, for $1 \le j \le t_i = \left \lceil \log |C_i| \right \rceil$.
Then there exists an index $h$, $1 \le h \le \left \lceil \log {|C_j| } 
\right \rceil $, such that the bit at the $h$-th position differs for $n_i(u)$
and $n_i(v)$. Without loss of generality assume that $h$-th bit in the binary
representation of $n_i(u)$  is $0$, and that of 
$n_i(v)$ is $1$.  Then by construction of  \rogers{$W_i^h$, 
$f_i^h (u) = [2,3]$ and  
$f_i^h(v) = [0,1]$}. It follows that  $(u,v) \notin E(W_i^h)$.
\end{case} 

\begin{case}[$u \in C_i$ and $v \in C_j$, where $i \ne j$] 
Let $h$, ($1 \le h \le \lceil \log \chi \rceil$),   be a position where
the binary representations of $i$ and $j$ differ.
Without loss of generality assume that the $h$-th bit in the binary
representation of \rogers{$i$ is $1$} and the $h$-th bit in the binary representation
of \rogers{$j$ is $0$}. Recalling the construction of the co-bipartite graph
$H_h$, we have $u \in A_h$ and $v \in B_h$ and therefore  $(u,v) \notin E(H_h)$.
Now, from equation \ref {cubicityofcobipartitegraph}, it follows that  there exists
$j$, $1 \le j \le 2b$, such
that in the unit interval graph $U_h^j$, $(u,v) \notin E(U_h^j)$.
\end{case} 

From Claim \ref {mainclaim} and Lemma \ref {Robertslem}, 
and noting that $t_i \le \left \lceil \log \alpha \right \rceil$ for
$0 \le i \le \chi -1$,
it immediately follows that
$\cubi(G) \le
 2 \left \lceil \log \chi(G)  \right \rceil \boxi(G) + 
   \chi(G)  \left \lceil \log \alpha(G)  \right \rceil $.

~~~~~~~

\noindent {\bf Tightness of Theorem \ref {mainthm} }:  
We will show that for every $\epsilon > 0$, there
exists a graph $G$ such that $\cubi(G) \le 
 2 \left \lceil \log \chi(G)  \right \rceil \boxi(G) + 
   \chi(G)  \left \lceil \log \alpha(G)  \right \rceil  \le \cubi(G) (1 + \epsilon)$. 
Let $k$ be a positive integer. 
Let $T_k$ be the complete  $k$-partite graph, where  each part
contains exactly  $\frac {n}{k}$ vertices. (Assume $n$ to be a multiple of 
$k$.)  The $\cubi(T_k) = k  \left \lceil \log \frac {n}{k} \right \rceil$. 
The upper bound of Theorem \ref {mainthm} for $T_k$ equals
$2  k  \left \lceil \log k  \right \rceil + k  \left \lceil \log \frac {n}{k} \right \rceil =  \cubi(T_k) ( 1 +  \frac {2 \left \lceil \log k \right \rceil}
   {\left \lceil \log \frac {n}{k} \right \rceil }  ) \le \cubi(T_k) (1 + \epsilon)$, provided we 
take  $n > k^{\frac {2 + \epsilon} {\epsilon} }$.  Thus there \rogers{exist}
graphs for which the upper bound given by Theorem \ref {mainthm} is arbitrarily
close to the true value of their cubicity.  

\bibliographystyle{plain}

\end {document}